\documentclass[11pt,a4paper]{article}

\usepackage[T1]{fontenc}
\usepackage[utf8]{inputenc}
\usepackage[sc,osf]{mathpazo}
\usepackage{amsmath,amssymb,amsthm,mathtools}
\usepackage[a4paper,top=30mm,bottom=28mm,left=28mm,right=28mm,headheight=14pt]{geometry}
\usepackage{microtype}
\usepackage[dvipsnames]{xcolor}
\usepackage{enumitem}
\usepackage{fancyhdr}
\usepackage{titlesec}
\usepackage{orcidlink}
\usepackage{hyperref}
\hypersetup{colorlinks=true,linkcolor=MidnightBlue,citecolor=OliveGreen,urlcolor=MidnightBlue,
  pdftitle={Stability of Caratheodory Outer Measures Generated by Spatial Gauges},
  pdfauthor={Rene Erlin Castillo},
  pdfsubject={Spatial gauges, Caratheodory outer measures, logarithmic stability, Baire-category localization},
  pdfkeywords={spatial gauge; metric outer measure; Hausdorff measure; stability; Baire category}}
\usepackage[capitalise,nameinlink]{cleveref}

\definecolor{accent}{RGB}{25,55,110}
\titleformat{\section}{\normalfont\Large\bfseries\color{accent}}{\thesection.}{0.6em}{}
\titleformat{\subsection}{\normalfont\large\bfseries\color{accent}}{\thesubsection.}{0.5em}{}
\titlespacing*{\section}{0pt}{2.2ex plus 1ex minus .2ex}{1.2ex plus .2ex}

\fancypagestyle{plain}{\fancyhf{}\fancyfoot[C]{\small\thepage}}

\setlist[enumerate]{leftmargin=2.2em,itemsep=2pt,topsep=4pt}
\setlist[itemize]{leftmargin=2.2em,itemsep=2pt,topsep=4pt}

\theoremstyle{plain}
\newtheorem{theorem}{Theorem}[section]
\newtheorem{lemma}[theorem]{Lemma}
\newtheorem{proposition}[theorem]{Proposition}
\newtheorem{corollary}[theorem]{Corollary}
\theoremstyle{definition}
\newtheorem{definition}[theorem]{Definition}
\newtheorem{example}[theorem]{Example}
\newtheorem{problem}[theorem]{Open problem}
\theoremstyle{remark}
\newtheorem{remark}[theorem]{Remark}
\crefname{problem}{Open problem}{Open problems}

\newcommand{\R}{\mathbb{R}}
\newcommand{\N}{\mathbb{N}}

\newcommand{\Hh}{\mathcal{H}}
\newcommand{\diam}{\operatorname{diam}}
\newcommand{\dist}{\operatorname{dist}}
\newcommand{\eps}{\varepsilon}
\newcommand{\cl}[1]{\overline{#1}}
\newcommand{\doi}[1]{\href{https://doi.org/#1}{doi:\nolinkurl{#1}}}

\begin{document}

\thispagestyle{plain}
\begin{center}
{\LARGE\bfseries Stability of Carathéodory Outer Measures}\\[0.6ex]
{\LARGE\bfseries Generated by Spatial Gauges}\\[1.2ex]
{\large Comparison, Logarithmic Perturbations, and Localization}\\[3.0ex]
{\large\scshape René Erlin Castillo}\,\orcidlink{0000-0003-1113-5827}\\[1.2ex]
{\small Departamento de Matemáticas, Universidad Nacional de Colombia,\\
Ciudad Universitaria, Bogotá D.C., Colombia\\[0.4ex]
\href{mailto:recastillo@unal.edu.co}{\texttt{recastillo@unal.edu.co}}
\quad$\cdot$\quad ORCID: \href{https://orcid.org/0000-0003-1113-5827}{\texttt{0000-0003-1113-5827}}\\}
{\large\scshape Héctor Camilo Chaparro}\,\orcidlink{0000-0002-0723-8199}\\[1.2ex]
{\small Programa de Matemáticas, Universidad de Cartagena,\\
	Cartagena de Indias, Colombia\\[0.4ex]
	\href{mailto:hchaparrog@unicartagena.edu.co}{\texttt{hchaparrog@unicartagena.edu.co}}
	\quad$\cdot$\quad ORCID: \href{https://orcid.org/0000-0002-0723-8199}{\texttt{0000-0002-0723-8199}}}
\end{center}

\vspace{1.5ex}
\begin{quote}
\small
\noindent\textbf{Abstract.}
We study outer measures generated by covering costs that depend on both the diameter and the
location of a covering set. After establishing the metric and Borel-regularity properties of this
construction, we investigate its stability under perturbations of the spatial gauge. For power gauges,
we prove a uniform logarithmic criterion: if the exponent perturbation multiplied by
$\log(1/r)$ converges uniformly to $\lambda$, then the induced measure is multiplied by
$e^{-\lambda}$. Pointwise asymptotic equivalence, however, is insufficient: on $[0,1]$ we construct
bounded gauges whose ratio is eventually equal to one at every fixed point, although their induced
measures differ by a prescribed factor. In the converse direction, continuity of the ratio in the
spatial variable on a Baire space yields equality of the measures on every compact subset of a
suitable dense $G_\delta$ set. The results identify a concrete obstruction to pointwise-to-measure
stability and a localization principle, while leaving open the corresponding question for all Borel
sets.


\end{quote}

\vspace{1ex}

\section{Introduction}

Hausdorff-type measures are constructed by assigning a cost to each member of a sufficiently fine
covering and then minimizing the total cost. This procedure is a form of Carathéodory's
\emph{Method~II} \cite{Caratheodory1914,Caratheodory1927}; its classical theory is developed in,
among other sources, \cite{Federer1969,Rogers1970,Mattila1995,Folland1999,Taylor2006,EoMOuter}.
In the standard Hausdorff construction, the cost depends only on the diameter of a covering set.
In several geometric settings, however, the cost also depends on the position of that set. Such
spatial dependence appears in multifractal constructions \cite{Olsen1995} and in the theory of
centred Hausdorff measures \cite{SaintRaymondTricot1988,Schechter2000,Schechter2002}.

The purpose of this paper is to examine how the outer measure changes when the spatial gauge is
perturbed, with particular emphasis on the distinction between uniform and pointwise asymptotic
control. Let $(X,d)$ be a metric space and let
$h\colon X\times(0,\infty)\to[0,\infty]$. For a set $A\subseteq X$ of positive finite diameter,
we assign the cost
\[
c_h(A)=\inf_{x\in A}h(x,\diam A),
\]
and assign cost zero to the empty set and to sets of diameter zero. The corresponding
$\delta$-contents are obtained by taking the infimum of the total cost over all countable
$\delta$-coverings; their limit as $\delta\downarrow0$ is denoted by $\mathcal H^h$.
The precise definitions are given in Section~\ref{sec:gauges}.

Some basic facts about this construction are needed before stability can be discussed. We prove
that each fixed-scale content is an outer measure and that the limiting set function is a metric
outer measure. We also prove Borel regularity, with an explicit $F_{\sigma\delta}$ hull for every
set. These statements provide the framework for the later arguments; the metric outer-measure
principle and Carathéodory's criterion themselves are classical, and we do not present them as
the principal source of novelty.

The main stability results are as follows.

\begin{enumerate}[label=\textup{(\roman*)}]
\item \textit{Uniform comparison and logarithmic perturbations.}
Uniform asymptotic equivalence of two spatial gauges implies equality of their limiting outer
measures. For gauges of the form
\[
h(x,r)=r^{s(x)},\qquad g(x,r)=r^{s(x)+a(x,r)},\quad 0<r<1,
\]
we obtain the sharper conclusion that
\[
\sup_{\substack{x\in X\\0<r\le\delta}}
\left|a(x,r)\log\frac1r-\lambda\right|\longrightarrow0
\quad\Longrightarrow\quad
\mathcal H^g=e^{-\lambda}\mathcal H^h.
\]
Thus the scale $1/\log(1/r)$ is the threshold for this uniform comparison principle.

\item \textit{Failure of pointwise stability.}
We construct bounded gauges on $[0,1]$ for which $g(x,r)/h(x,r)\to1$ at every fixed point—in
fact, the ratio is identically one for all sufficiently small $r$, with the threshold depending
on the point—while
\[
\mathcal H^g(E)=\kappa\,\mathcal H^h(E)
\qquad\text{for every }E\subseteq[0,1],
\]
where $\kappa\in(0,1)$ is arbitrary. This shows that pointwise asymptotic equivalence alone
cannot control the induced outer measure.

\item \textit{Baire-category localization.}
If $g(\cdot,r)/h(\cdot,r)$ is continuous for every $r$ and converges pointwise to one on a
Baire space, then for each $\eta>0$ there is a dense open set on which the convergence is
locally uniform to within $\eta$. Intersecting these sets gives a dense $G_\delta$ set $G$.
We then prove equality $\mathcal H^g(K)=\mathcal H^h(K)$ for every compact $K\subseteq G$.
More generally, a two-sided estimate holds on compact subsets of the corresponding
$\eta$-localization set.
\end{enumerate}

The last result is deliberately a localization theorem, not a theorem for all Borel sets.
The exceptional set $X\setminus G$ is meagre, but it need not be negligible for the outer
measure under consideration. Consequently, category-smallness alone does not settle the
pointwise-to-measure problem on arbitrary Borel sets. We formulate this remaining issue as an
open problem in Section~\ref{sec:local}.

The paper is organized as follows. Section~\ref{sec:gauges} introduces the spatial gauges and
covering contents. Section~\ref{sec:metric} establishes the metric outer-measure property.
Section~\ref{sec:comparison} develops comparison principles, and Section~\ref{sec:power}
contains the logarithmic criterion. Section~\ref{sec:regularity} discusses measurability and
regularity. The counterexample and the Baire-category localization theorem are proved in
Section~\ref{sec:local}. Section~\ref{sec:literature} places these results in the context of
spatially dependent and centred Hausdorff measures.

\medskip
\noindent\textbf{Conventions.}
We work in $[0,\infty]$ with the usual order and the conventions
$a+\infty=\infty$, $a\cdot\infty=\infty$ for $a>0$, $0\cdot\infty=0$, and
$\inf\emptyset=\infty$. For $A\subseteq X$ we put $\diam\emptyset=0$ and
$\diam A=\sup\{d(x,y):x,y\in A\}$ otherwise. For nonempty $E,F\subseteq X$,
$\dist(E,F)=\inf\{d(x,y):x\in E,\ y\in F\}$. The Lebesgue outer measure on
$\mathbb R$ is denoted by $\lambda^*$; we use its standard properties, including
$\lambda^*(I)=|I|$ for every interval $I$, as developed for example in
\cite[Chapter~1]{Folland1999}.

\section{Spatial gauges and covering constructions}\label{sec:gauges}

Throughout, $(X,d)$ is a metric space.

\begin{definition}[Spatial gauge]\label{def:gauge}
A \emph{spatial gauge} is a function $h\colon X\times(0,\infty)\to[0,\infty]$.
No continuity, measurability or monotonicity is assumed unless explicitly stated.
\end{definition}

\begin{definition}[Covering cost and the measures $\Hh^h_\delta$, $\Hh^h$]\label{def:cost}
For $A\subseteq X$ define
\[
c_h(A)=
\begin{cases}
\displaystyle\inf_{x\in A} h(x,\diam A), & A\neq\emptyset,\ 0<\diam A<\infty,\\[1.2ex]
0, & A=\emptyset \text{ or } \diam A=0 .
\end{cases}
\]
(Sets of infinite diameter never occur in what follows.) For $\delta>0$, a sequence $(A_j)_{j\ge1}$ of
subsets of $X$ is a \emph{$\delta$-covering} of $E\subseteq X$ if $E\subseteq\bigcup_j A_j$ and
$\diam A_j\le\delta$ for all $j$. Set
\[
\Hh^h_\delta(E)=\inf\Bigl\{\sum_{j=1}^{\infty}c_h(A_j)\ :\ (A_j)_{j\ge1}\ \text{is a $\delta$-covering of } E\Bigr\},
\qquad
\Hh^h(E)=\lim_{\delta\downarrow0}\Hh^h_\delta(E).
\]
\end{definition}

Finite coverings are included, since a finite family may be completed by empty sets, which have
diameter and cost $0$. If $E$ admits no countable $\delta$-covering (which can happen only if $X$ is
not separable), then $\Hh^h_\delta(E)=\inf\emptyset=\infty$.

\begin{lemma}[Existence of the limit]\label{lem:limit}
For every $E\subseteq X$ the function $\delta\mapsto\Hh^h_\delta(E)$ is nonincreasing on $(0,\infty)$.
Consequently the limit $\Hh^h(E)$ exists in $[0,\infty]$ and
$\Hh^h(E)=\sup_{\delta>0}\Hh^h_\delta(E)$.
\end{lemma}

\begin{proof}
Let $0<\delta_1<\delta_2$. Every $\delta_1$-covering of $E$ is a $\delta_2$-covering of $E$, so the infimum
defining $\Hh^h_{\delta_2}(E)$ is taken over a larger class than the one defining
$\Hh^h_{\delta_1}(E)$. Hence $\Hh^h_{\delta_2}(E)\le\Hh^h_{\delta_1}(E)$. A nonincreasing function of
$\delta$ with values in $[0,\infty]$ has a limit as $\delta\downarrow0$, equal to its supremum.
\end{proof}

\begin{remark}[Sets of diameter zero]\label{rem:diamzero}
The convention $c_h(A)=0$ when $\diam A=0$ is the usual Hausdorff-type convention for gauges satisfying
$h(x,r)\to0$ as $r\downarrow0$; it is imposed here explicitly so that no value $h(x,0)$ is needed.
A set of diameter zero is either empty or a singleton. Consequently every countable set
$E=\{x_1,x_2,\dots\}$ satisfies $\Hh^h(E)=0$: for every $\delta>0$ the sequence of singletons
$(\{x_j\})_{j\ge1}$ is a $\delta$-covering of $E$ of total cost $0$, so $\Hh^h_\delta(E)=0$.
\end{remark}

The next two lemmas are used repeatedly.

\begin{lemma}[Locality in the radius]\label{lem:locality}
Let $h,g$ be spatial gauges and $r_0>0$ with $h(x,r)=g(x,r)$ for all $x\in X$ and $0<r<r_0$. Then
$\Hh^h_\delta=\Hh^g_\delta$ for $0<\delta<r_0$, and hence $\Hh^h=\Hh^g$. Moreover, $\Hh^{\kappa h}=\kappa\,\Hh^h$
for every constant $\kappa\in(0,\infty)$.
\end{lemma}

\begin{proof}
If $0<\delta<r_0$ and $\diam A\le\delta$, then either $\diam A=0$ and $c_h(A)=0=c_g(A)$, or
$0<\diam A<r_0$ and $c_h(A)=c_g(A)$ because the two infima involve the same values. Hence the two
infima defining $\Hh^h_\delta(E)$ and $\Hh^g_\delta(E)$ coincide, and \cref{lem:limit} gives
$\Hh^h(E)=\lim_{\delta\downarrow0}\Hh^h_\delta(E)=\lim_{\delta\downarrow0}\Hh^g_\delta(E)=\Hh^g(E)$.
For the last assertion, $c_{\kappa h}(A)=\kappa c_h(A)$ for every $A$ (an infimum commutes with
multiplication by $\kappa>0$), hence $\sum_j c_{\kappa h}(A_j)=\kappa\sum_j c_h(A_j)$ for every covering;
taking infima and then limits gives $\Hh^{\kappa h}_\delta=\kappa\Hh^h_\delta$ and
$\Hh^{\kappa h}=\kappa\Hh^h$.
\end{proof}

In particular, a gauge need only be specified for small radii; we shall use this freely.

\begin{lemma}[Closed coverings suffice]\label{lem:closed}
For every $A\subseteq X$ one has $\diam\cl A=\diam A$ and $c_h(\cl A)\le c_h(A)$. Consequently, for every
$E\subseteq X$ and $\delta>0$,
\[
\Hh^h_\delta(E)=\inf\Bigl\{\sum_{j}c_h(F_j)\ :\ (F_j)_{j\ge1}\ \text{is a $\delta$-covering of $E$ by closed sets}\Bigr\}.
\]
\end{lemma}

\begin{proof}
If $x,y\in\cl A$, choose $x_n,y_n\in A$ with $x_n\to x$, $y_n\to y$; then
$d(x,y)=\lim d(x_n,y_n)\le\diam A$. Thus $\diam\cl A\le\diam A$, and the reverse inequality is obvious.
If $\diam A=0$, then $A$ is empty or a singleton, hence closed, and $\cl A=A$. If $r=\diam A>0$, then
$\diam\cl A=r$ and, since $A\subseteq\cl A$,
\[
c_h(\cl A)=\inf_{x\in\cl A}h(x,r)\le\inf_{x\in A}h(x,r)=c_h(A).
\]
Let $I$ denote the infimum over closed $\delta$-coverings. Every closed $\delta$-covering is a
$\delta$-covering, so $\Hh^h_\delta(E)\le I$. Conversely, if $(A_j)$ is any $\delta$-covering of $E$, then
$(\cl{A_j})$ is a closed $\delta$-covering of $E$ with $\sum_j c_h(\cl{A_j})\le\sum_j c_h(A_j)$, so
$I\le\sum_jc_h(A_j)$; taking the infimum over $(A_j)$ gives $I\le\Hh^h_\delta(E)$.
\end{proof}

\begin{remark}
The inequality $c_h(\cl A)\le c_h(A)$ is a feature of the present construction, in which the location
parameter is taken from the \emph{covering set}. It fails for centred constructions, where the centres
must lie in the set being measured \cite{SaintRaymondTricot1988,Schechter2000}; this is the source of
the difficulty of Borel regularity for centred Hausdorff measures (see \cref{sec:literature}).
\end{remark}

\section{The metric outer-measure theorem}\label{sec:metric}

Recall that an \emph{outer measure} on $X$ is a map $\mu\colon\mathcal P(X)\to[0,\infty]$ with
$\mu(\emptyset)=0$, which is monotone and countably subadditive; a set $M$ is
\emph{$\mu$-measurable} (in the sense of Carathéodory) if $\mu(T)=\mu(T\cap M)+\mu(T\setminus M)$ for every
$T\subseteq X$. An outer measure is \emph{metric} if $\mu(E\cup F)=\mu(E)+\mu(F)$ whenever
$\dist(E,F)>0$.

\begin{theorem}[Metric Carathéodory construction]\label{thm:metric}
For every spatial gauge $h$, each $\Hh^h_\delta$ $(\delta>0)$ is an outer measure, and $\Hh^h$ is a metric
outer measure on $X$. In particular, every Borel subset of $X$ is $\Hh^h$-measurable in the sense of
Carathéodory, and the restriction of $\Hh^h$ to the Borel $\sigma$-algebra is a countably additive
measure.
\end{theorem}

\begin{proof}
\emph{Step 1: $\Hh^h_\delta(\emptyset)=0$ and $\Hh^h(\emptyset)=0$.} The sequence $A_j=\emptyset$ $(j\ge1)$ is a
$\delta$-covering of $\emptyset$ of cost $0$; hence $\Hh^h_\delta(\emptyset)=0$ for every $\delta$, and
$\Hh^h(\emptyset)=0$.

\emph{Step 2: monotonicity.} If $E\subseteq F$, every $\delta$-covering of $F$ is a $\delta$-covering of $E$,
so $\Hh^h_\delta(E)\le\Hh^h_\delta(F)$ for each $\delta>0$. Letting $\delta\downarrow0$ gives
$\Hh^h(E)\le\Hh^h(F)$.

\emph{Step 3: countable subadditivity of $\Hh^h_\delta$.} Fix $\delta>0$ and sets $E_k\subseteq X$,
$k\in\N$. If $\sum_k\Hh^h_\delta(E_k)=\infty$ there is nothing to prove. Otherwise every
$\Hh^h_\delta(E_k)$ is finite; given $\eps>0$, by definition of the infimum we can choose, for each $k$,
a $\delta$-covering $(A_{k,j})_{j\ge1}$ of $E_k$ such that
\[
\sum_{j\ge1}c_h(A_{k,j})\le\Hh^h_\delta(E_k)+\eps\,2^{-k}.
\]
Enumerate the countable family $\{A_{k,j}:k,j\ge1\}$ as a single sequence (for instance via a bijection
$\N\to\N\times\N$). It covers $\bigcup_kE_k$, each of its members has diameter at most $\delta$, and,
because the terms are nonnegative, its total cost does not depend on the enumeration. Therefore
\[
\Hh^h_\delta\Bigl(\bigcup_kE_k\Bigr)\le\sum_{k}\sum_{j}c_h(A_{k,j})
\le\sum_k\Hh^h_\delta(E_k)+\eps\sum_k2^{-k}=\sum_k\Hh^h_\delta(E_k)+\eps .
\]
Letting $\eps\downarrow0$ yields $\Hh^h_\delta(\bigcup_kE_k)\le\sum_k\Hh^h_\delta(E_k)$. Together with
Steps~1 and~2, this shows that $\Hh^h_\delta$ is an outer measure.

\emph{Step 4: countable subadditivity of $\Hh^h$.} By \cref{lem:limit}, $\Hh^h_\delta(E_k)\le\Hh^h(E_k)$ for
every $k$ and every $\delta>0$. Hence, by Step~3,
\[
\Hh^h_\delta\Bigl(\bigcup_kE_k\Bigr)\le\sum_k\Hh^h_\delta(E_k)\le\sum_k\Hh^h(E_k)
\qquad(\delta>0).
\]
The right-hand side does not depend on $\delta$; letting $\delta\downarrow0$ on the left gives
$\Hh^h(\bigcup_kE_k)\le\sum_k\Hh^h(E_k)$. Thus $\Hh^h$ is an outer measure.

\emph{Step 5: the metric property.} Let $E,F\subseteq X$ with $a:=\dist(E,F)>0$ and fix
$0<\delta<a$. If $E\cup F$ has no $\delta$-covering, then $\Hh^h_\delta(E\cup F)=\infty$ and the inequality
\eqref{eq:metricdelta} below is trivial. Otherwise, let $(A_j)$ be a $\delta$-covering of $E\cup F$. No
$A_j$ meets both $E$ and $F$: if $x\in A_j\cap E$ and $y\in A_j\cap F$, then
$a\le d(x,y)\le\diam A_j\le\delta<a$, a contradiction. Put
\[
I_E=\{j:A_j\cap E\neq\emptyset\},\qquad I_F=\{j:A_j\cap F\neq\emptyset\}.
\]
These index sets are disjoint; the subfamily $(A_j)_{j\in I_E}$ covers $E$ and the subfamily
$(A_j)_{j\in I_F}$ covers $F$, and both are $\delta$-coverings (after completion by empty sets). Since all
costs are nonnegative,
\[
\sum_jc_h(A_j)\ \ge\ \sum_{j\in I_E}c_h(A_j)+\sum_{j\in I_F}c_h(A_j)\ \ge\ \Hh^h_\delta(E)+\Hh^h_\delta(F).
\]
Taking the infimum over all $\delta$-coverings of $E\cup F$ yields
\begin{equation}\label{eq:metricdelta}
\Hh^h_\delta(E\cup F)\ \ge\ \Hh^h_\delta(E)+\Hh^h_\delta(F)\qquad(0<\delta<a).
\end{equation}
The reverse inequality holds by Step~3, so equality holds in \eqref{eq:metricdelta} for every
$\delta<a$. Letting $\delta\downarrow0$ (both sides are nondecreasing as $\delta$ decreases, so the limit of
the sum is the sum of the limits in $[0,\infty]$) we obtain $\Hh^h(E\cup F)=\Hh^h(E)+\Hh^h(F)$.

\emph{Step 6: Borel sets.} By Steps~1--5, $\Hh^h$ is a metric outer measure. By
\cref{lem:caratheodory} below, every closed set is $\Hh^h$-measurable. By Carathéodory's theorem
\cite{Caratheodory1914}, \cite[Theorem~1.11]{Folland1999}, the $\Hh^h$-measurable sets form a
$\sigma$-algebra on which $\Hh^h$ is countably additive. This $\sigma$-algebra contains all closed sets,
hence all open sets, hence the $\sigma$-algebra they generate, which is the Borel $\sigma$-algebra.
\end{proof}

For completeness we include the proof of Carathéodory's criterion in the form we need.

\begin{lemma}[Carathéodory's criterion]\label{lem:caratheodory}
Let $\mu$ be a metric outer measure on $X$. Then every closed set $C\subseteq X$ is $\mu$-measurable.
\end{lemma}

\begin{proof}
Let $T\subseteq X$. By subadditivity $\mu(T)\le\mu(T\cap C)+\mu(T\setminus C)$, so it suffices to prove
\begin{equation}\label{eq:caraineq}
\mu(T\cap C)+\mu(T\setminus C)\le\mu(T),
\end{equation}
and we may assume $\mu(T)<\infty$. If $C=\emptyset$, \eqref{eq:caraineq} is trivial, so let
$C\ne\emptyset$ and write $\rho(x)=\dist(\{x\},C)$. For $n\in\N$ put
\[
T_n=\{x\in T\setminus C:\ \rho(x)\ge 1/n\},\qquad
R_k=\Bigl\{x\in T\setminus C:\ \tfrac1{k+1}\le\rho(x)<\tfrac1k\Bigr\}\quad(k\in\N).
\]
\emph{(a)} Since $\dist(T\cap C,T_n)\ge1/n>0$ (if both are nonempty), the metric property and monotonicity give
$\mu(T\cap C)+\mu(T_n)=\mu((T\cap C)\cup T_n)\le\mu(T)$ for every $n$.

\emph{(b)} Since $C$ is closed, $\rho(x)>0$ for every $x\notin C$; hence
$T\setminus C=T_n\cup\bigcup_{k\ge n}R_k$ for every $n$.

\emph{(c)} If $k\ge j+2$, then $\dist(R_j,R_k)>0$: for $x\in R_j$ and $y\in R_k$, the triangle
inequality gives $\rho(x)\le d(x,y)+\rho(y)$, so
$d(x,y)\ge\rho(x)-\rho(y)>\frac1{j+1}-\frac1k\ge\frac1{j+1}-\frac1{j+2}>0$. Consequently, for every
$N$, the sets $R_2,R_4,\dots,R_{2N}$ are pairwise positively separated, and so are
$R_1,R_3,\dots,R_{2N-1}$. Since the distance from a finite union to a set is the minimum of the
distances, an induction on $N$ using the metric property gives
\[
\sum_{i=1}^N\mu(R_{2i})=\mu\Bigl(\bigcup_{i=1}^NR_{2i}\Bigr)\le\mu(T),\qquad
\sum_{i=1}^N\mu(R_{2i-1})=\mu\Bigl(\bigcup_{i=1}^NR_{2i-1}\Bigr)\le\mu(T).
\]
Hence $\sum_{k\ge1}\mu(R_k)\le2\mu(T)<\infty$.

\emph{(d)} By (b) and subadditivity, $\mu(T\setminus C)\le\mu(T_n)+\sum_{k\ge n}\mu(R_k)$ for every $n$. The
tail of a convergent series tends to $0$, and $\mu(T_n)$ is nondecreasing in $n$ and bounded by
$\mu(T\setminus C)$; therefore $\mu(T\setminus C)=\lim_n\mu(T_n)$. Letting $n\to\infty$ in (a) gives
\eqref{eq:caraineq}.
\end{proof}

\begin{remark}[What \cref{thm:metric} does not say]\label{rem:notsay}
The theorem does not require $h$ to be measurable in $x$, continuous, or monotone in $r$. It does not
establish local finiteness or any inner regularity. Borel regularity, which does not follow from the
metric property alone, is proved separately in \cref{thm:borelregular}; outer regularity by open sets
fails in general (\cref{ex:noopen}).
\end{remark}

\section{Comparison of spatial gauges}\label{sec:comparison}

We first isolate the elementary comparison mechanism.

\begin{lemma}[Comparison at a fixed scale]\label{lem:fixedscale}
Let $h,g$ be spatial gauges, $\delta>0$ and $0<m\le M<\infty$ such that
\[
m\,h(x,r)\le g(x,r)\le M\,h(x,r)\qquad\text{for all }x\in X,\ 0<r\le\delta.
\]
Then $m\,c_h(A)\le c_g(A)\le M\,c_h(A)$ whenever $\diam A\le\delta$, and
$m\,\Hh^h_\delta(E)\le\Hh^g_\delta(E)\le M\,\Hh^h_\delta(E)$ for every $E\subseteq X$.
\end{lemma}

\begin{proof}
If $\diam A=0$ all three costs vanish. Let $0<r=\diam A\le\delta$. For every $x\in A$,
$g(x,r)\ge m\,h(x,r)\ge m\,c_h(A)$; taking the infimum over $x\in A$ gives $c_g(A)\ge m\,c_h(A)$.
Also $c_g(A)\le g(x,r)\le M\,h(x,r)$ for every $x\in A$, and taking the infimum over $x\in A$ gives
$c_g(A)\le M\inf_{x\in A}h(x,r)=M\,c_h(A)$.

Now let $E\subseteq X$. If $E$ has no $\delta$-covering, all three quantities are $\infty$ (recall $m>0$).
Otherwise, for every $\delta$-covering $(A_j)$ of $E$,
\[
\sum_jc_g(A_j)\ge m\sum_jc_h(A_j)\ge m\,\Hh^h_\delta(E),
\qquad\text{and}\qquad
\Hh^g_\delta(E)\le\sum_jc_g(A_j)\le M\sum_jc_h(A_j).
\]
Taking the infimum over all $\delta$-coverings in the first chain gives $\Hh^g_\delta(E)\ge m\Hh^h_\delta(E)$;
in the second chain it gives $\Hh^g_\delta(E)\le M\Hh^h_\delta(E)$, because
$\inf(M\,t)=M\inf t$ for $M\in(0,\infty)$.
\end{proof}

\begin{lemma}[Passage to the limit]\label{lem:passage}
Let $h,g$ be spatial gauges, $\delta_0>0$ and $m,M\colon(0,\delta_0]\to(0,\infty)$ with
$m(\delta)\to1$ and $M(\delta)\to1$ as $\delta\downarrow0$, such that
$m(\delta)\Hh^h_\delta(E)\le\Hh^g_\delta(E)\le M(\delta)\Hh^h_\delta(E)$ for all $E\subseteq X$ and
$0<\delta\le\delta_0$. Then $\Hh^g=\Hh^h$.
\end{lemma}

\begin{proof}
Fix $E$. If $\Hh^h(E)<\infty$, then $\Hh^h_\delta(E)\to\Hh^h(E)$ in $[0,\infty)$, so both
$m(\delta)\Hh^h_\delta(E)$ and $M(\delta)\Hh^h_\delta(E)$ converge to $\Hh^h(E)$, and the squeeze gives
$\Hh^g(E)=\Hh^h(E)$. If $\Hh^h(E)=\infty$, then $\Hh^h_\delta(E)\to\infty$; since $m(\delta)\ge\frac12$ for
small $\delta$, the lower bound forces $\Hh^g_\delta(E)\to\infty$, i.e.\ $\Hh^g(E)=\infty$. (Symmetrically,
$\Hh^g(E)=\infty$ forces $\Hh^h(E)=\infty$ through the upper bound, so the two cases are exhaustive.)
\end{proof}

\begin{theorem}[Uniform asymptotic equivalence]\label{thm:uniform}
Let $h,g$ be spatial gauges. Suppose there are $r_0>0$ and a function $\eps\colon(0,r_0)\to[0,1)$ with
$\lim_{r\downarrow0}\eps(r)=0$ such that
\[
(1-\eps(r))\,h(x,r)\le g(x,r)\le(1+\eps(r))\,h(x,r)\qquad\text{for all }x\in X,\ 0<r<r_0.
\]
Then $\Hh^h(E)=\Hh^g(E)$ for every $E\subseteq X$.
\end{theorem}

\begin{proof}
For $0<\delta<r_0$ put $\eta(\delta)=\sup_{0<r\le\delta}\eps(r)\in[0,1]$. The function $\eta$ is
nondecreasing, and $\eta(\delta)\to0$ as $\delta\downarrow0$: given $\theta>0$ there is $\rho>0$ with
$\eps(r)<\theta$ for $0<r<\rho$, hence $\eta(\delta)\le\theta$ for $\delta<\rho$. Choose $\delta_0<r_0$ with
$\eta(\delta_0)\le\frac12$. For $0<\delta\le\delta_0$ and $0<r\le\delta$ we have
$\eps(r)\le\eta(\delta)$, hence
\[
(1-\eta(\delta))\,h(x,r)\le g(x,r)\le(1+\eta(\delta))\,h(x,r)\qquad(x\in X).
\]
\Cref{lem:fixedscale} with $m=1-\eta(\delta)\ge\frac12$ and $M=1+\eta(\delta)$ gives
$(1-\eta(\delta))\Hh^h_\delta\le\Hh^g_\delta\le(1+\eta(\delta))\Hh^h_\delta$, and \cref{lem:passage} yields
$\Hh^g=\Hh^h$, including the case in which the common value is infinite.
\end{proof}

\begin{corollary}[Two-sided comparability]\label{cor:comparable}
Suppose there are constants $0<m\le M<\infty$ and $r_0>0$ such that
$m\,h(x,r)\le g(x,r)\le M\,h(x,r)$ for all $x\in X$ and $0<r<r_0$. Then
\[
m\,\Hh^h(E)\le\Hh^g(E)\le M\,\Hh^h(E)\qquad\text{for every }E\subseteq X.
\]
\end{corollary}

\begin{proof}
For $0<\delta<r_0$, \cref{lem:fixedscale} gives $m\Hh^h_\delta(E)\le\Hh^g_\delta(E)\le M\Hh^h_\delta(E)$.
Since $\Hh^h_\delta(E)\uparrow\Hh^h(E)$ and $\Hh^g_\delta(E)\uparrow\Hh^g(E)$ as $\delta\downarrow0$, and
multiplication by a positive constant is continuous and nondecreasing on $[0,\infty]$, letting
$\delta\downarrow0$ gives the claim.
\end{proof}

\section{Power gauges and the logarithmic threshold}\label{sec:power}

We now specialize to power-type gauges. Let $s\colon X\to(0,\infty)$ and, for $0<r<1$, put
$h(x,r)=r^{s(x)}$. Let $a\colon X\times(0,1)\to\R$ and put $g(x,r)=r^{s(x)+a(x,r)}$ for $0<r<1$. For
$r\ge1$ both gauges may be defined arbitrarily: by \cref{lem:locality} this does not affect $\Hh^h$ or
$\Hh^g$. (The condition $s(x)+a(x,r)>0$, natural if one wants $g(x,r)\to0$, is not needed below.)
The basic identity is
\begin{equation}\label{eq:ratio}
\frac{g(x,r)}{h(x,r)}=r^{a(x,r)}=\exp\bigl(a(x,r)\log r\bigr)=\exp\Bigl(-a(x,r)\log\frac1r\Bigr),
\qquad 0<r<1 .
\end{equation}

\begin{theorem}[Logarithmic comparison]\label{thm:log}
Let $\lambda\in\R$ and assume
\[
\theta(\delta):=\sup\Bigl\{\Bigl|a(x,r)\log\frac1r-\lambda\Bigr| : x\in X,\ 0<r\le\delta\Bigr\}\longrightarrow0
\qquad(\delta\downarrow0).
\]
Then $\Hh^g(E)=e^{-\lambda}\,\Hh^h(E)$ for every $E\subseteq X$.
\end{theorem}

\begin{proof}
Let $\tilde h=e^{-\lambda}h$. By \eqref{eq:ratio},
$g(x,r)/\tilde h(x,r)=\exp\bigl(-(a(x,r)\log\frac1r-\lambda)\bigr)$ for $0<r<1$. Choose
$\delta_0\in(0,1)$ with $\theta(\delta_0)<\infty$. For $0<\delta\le\delta_0$, $x\in X$ and $0<r\le\delta$,
the exponent has absolute value at most $\theta(\delta)$, so
\[
e^{-\theta(\delta)}\,\tilde h(x,r)\le g(x,r)\le e^{\theta(\delta)}\,\tilde h(x,r).
\]
(Here $\tilde h(x,r)\in(0,\infty)$.) \Cref{lem:fixedscale} gives
$e^{-\theta(\delta)}\Hh^{\tilde h}_\delta\le\Hh^g_\delta\le e^{\theta(\delta)}\Hh^{\tilde h}_\delta$, and since
$e^{\pm\theta(\delta)}\to1$, \cref{lem:passage} yields $\Hh^g=\Hh^{\tilde h}$. Finally
$\Hh^{\tilde h}=e^{-\lambda}\Hh^h$ by \cref{lem:locality}.
\end{proof}

\begin{corollary}[Logarithmic stability criterion]\label{cor:logstab}
If
\[
\lim_{\delta\downarrow0}\ \sup_{x\in X,\ 0<r\le\delta}|a(x,r)|\log\frac1r=0,
\]
then $\Hh^g(E)=\Hh^h(E)$ for every $E\subseteq X$.
\end{corollary}

\begin{proof}
This is \cref{thm:log} with $\lambda=0$. Alternatively, put
$\omega(\delta)=\sup_{x,\,0<r\le\delta}|a(x,r)|\log\frac1r$; then \eqref{eq:ratio} gives
$e^{-\omega(\delta)}\le g/h\le e^{\omega(\delta)}$ for $0<r\le\delta$, and one concludes with
\cref{lem:fixedscale,lem:passage}.
\end{proof}

To evaluate the measures in the example below we need the classical identification of the
one-dimensional measure on the line.

\begin{lemma}\label{lem:lebesgue}
Let $X\subseteq\R$ be an interval with the Euclidean metric and $h(x,r)=r$. Then $\Hh^h_\delta(E)=\lambda^*(E)$
for every $E\subseteq X$ and every $\delta>0$; in particular $\Hh^h(E)=\lambda^*(E)$, and
$\Hh^h([0,1])=1$ when $X=[0,1]$.
\end{lemma}

\begin{proof}
Here $c_h(A)=\diam A$ for every $A\subseteq X$.

\emph{Lower bound.} Let $(A_j)$ be a $\delta$-covering of $E$. For $A_j\ne\emptyset$ let
$I_j=[\inf A_j,\sup A_j]$, a closed interval with $A_j\subseteq I_j$ and $|I_j|=\diam A_j$; for
$A_j=\emptyset$ let $I_j=\emptyset$. Then $E\subseteq\bigcup_jI_j$ and, by monotonicity and countable
subadditivity of $\lambda^*$,
$\lambda^*(E)\le\sum_j\lambda^*(I_j)=\sum_j\diam A_j=\sum_jc_h(A_j)$. Taking the infimum gives
$\lambda^*(E)\le\Hh^h_\delta(E)$.

\emph{Upper bound.} We may assume $\lambda^*(E)<\infty$. Let $\eps>0$ and choose open intervals $J_k$
with $E\subseteq\bigcup_kJ_k$ and $\sum_k|J_k|\le\lambda^*(E)+\eps$; we may discard empty $J_k$, and
since $\sum_k|J_k|<\infty$ each $J_k$ is bounded. Divide the closure of $J_k$ into $N_k$ consecutive closed
intervals $P_{k,1},\dots,P_{k,N_k}$ of equal length $|J_k|/N_k\le\delta$. The sets $P_{k,i}\cap X$ are
subsets of $X$ with $\diam(P_{k,i}\cap X)\le|P_{k,i}|\le\delta$, they cover $E$, and
\[
\sum_{k,i}c_h(P_{k,i}\cap X)\le\sum_{k}\sum_{i=1}^{N_k}|P_{k,i}|=\sum_k|J_k|\le\lambda^*(E)+\eps .
\]
Hence $\Hh^h_\delta(E)\le\lambda^*(E)+\eps$ for every $\eps>0$. Combining both bounds gives
$\Hh^h_\delta(E)=\lambda^*(E)$ for every $\delta$, and the remaining assertions follow from
$\lambda^*([0,1])=1$ \cite[Chapter~1]{Folland1999}.
\end{proof}

\begin{example}[A perturbation at the logarithmic scale]\label{ex:log}
Let $X=[0,1]$ with the Euclidean metric, $h(x,r)=r$, and fix $a\in\R$. Put $r_a=\min\{1,e^{a}\}$ and
define, for $0<r<r_a$,
\[
g(x,r)=r^{\,1+a/\log(1/r)} .
\]
(For $0<r<r_a$ we have $\log(1/r)>\max\{0,-a\}$, so the exponent $1+a/\log(1/r)$ is positive; for
$r\ge r_a$ the gauge $g$ may be defined arbitrarily, by \cref{lem:locality}.) Since
$r^{a/\log(1/r)}=\exp\bigl(\frac{a}{\log(1/r)}\log r\bigr)=\exp(-a)$, we have
\[
g(x,r)=e^{-a}\,r\qquad(0<r<r_a).
\]
By \cref{lem:locality,lem:lebesgue}, $\Hh^g(E)=e^{-a}\Hh^h(E)=e^{-a}\lambda^*(E)$ for every $E\subseteq[0,1]$.
In particular,
\[
\Hh^h([0,1])=1,\qquad\Hh^g([0,1])=e^{-a}.
\]
Here the exponent perturbation $a(x,r)=a/\log(1/r)$ tends to zero, but
$a(x,r)\log(1/r)\equiv a$; this agrees with \cref{thm:log} for $\lambda=a$. Thus a perturbation of the
exponent that tends to zero need not preserve the measure when it is of order $1/\log(1/r)$.
\end{example}

\begin{remark}[Sharpness]\label{rem:sharp}
\Cref{cor:logstab} and \cref{ex:log} identify the sharp scale for the uniform comparison method:
perturbations of order $o(1/\log(1/r))$, uniformly in $x$, preserve the measure, whereas a perturbation
equal to $\lambda/\log(1/r)$ with $\lambda\ne0$ multiplies it by $e^{-\lambda}\neq1$
(\cref{thm:log}). This is a sharpness statement for uniform criteria; it is not a classification of all
possible nonuniform perturbations (see \cref{sec:local}).
\end{remark}

\section{Regularity: measurability, Borel regularity and open sets}\label{sec:regularity}

\begin{definition}\label{def:borelreg}
An outer measure $\mu$ on $X$ is \emph{Borel regular} if every Borel set is $\mu$-measurable and for every
$E\subseteq X$ there exists a Borel set $B\supseteq E$ with $\mu(B)=\mu(E)$.
\end{definition}

The metric property (\cref{thm:metric}) yields the first requirement but, by itself, gives no Borel hull:
the proof of \cref{thm:metric} constructs no approximation from above. The hull is provided by the closure
lemma (\cref{lem:closed}).

\begin{theorem}[Borel regularity]\label{thm:borelregular}
For every spatial gauge $h$ the outer measure $\Hh^h$ is Borel regular. More precisely, for every
$E\subseteq X$ there is a set $B\supseteq E$ which is a countable intersection of countable unions of
closed sets (an $F_{\sigma\delta}$ set) such that $\Hh^h(B)=\Hh^h(E)$.
\end{theorem}

\begin{proof}
By \cref{thm:metric} every Borel set is $\Hh^h$-measurable. Let $E\subseteq X$. If $\Hh^h(E)=\infty$,
take $B=X$ (which is closed), since $\Hh^h(X)\ge\Hh^h(E)=\infty$. Assume $\Hh^h(E)<\infty$; then
$\Hh^h_\delta(E)\le\Hh^h(E)<\infty$ for every $\delta$ (\cref{lem:limit}). For each $n\in\N$, by
\cref{lem:closed} there is a $1/n$-covering $(F_{n,j})_{j\ge1}$ of $E$ by closed sets with
\[
\sum_jc_h(F_{n,j})\le\Hh^h_{1/n}(E)+\tfrac1n\le\Hh^h(E)+\tfrac1n .
\]
Put $B_n=\bigcup_jF_{n,j}$ and $B=\bigcap_nB_n$. Then $B$ is $F_{\sigma\delta}$, hence Borel, and
$E\subseteq B$ because $E\subseteq B_n$ for every $n$. Since $B\subseteq B_n$, the family $(F_{n,j})_j$ is
a $1/n$-covering of $B$, whence
\[
\Hh^h_{1/n}(B)\le\sum_jc_h(F_{n,j})\le\Hh^h(E)+\tfrac1n\qquad(n\in\N).
\]
Fix $m\in\N$. For $n\ge m$, \cref{lem:limit} gives $\Hh^h_{1/m}(B)\le\Hh^h_{1/n}(B)\le\Hh^h(E)+\frac1n$;
letting $n\to\infty$ we get $\Hh^h_{1/m}(B)\le\Hh^h(E)$, and letting $m\to\infty$ we get
$\Hh^h(B)\le\Hh^h(E)$. The reverse inequality is monotonicity.
\end{proof}

As usual, Borel regularity yields continuity along increasing sequences of arbitrary sets.

\begin{corollary}\label{cor:increasing}
If $E_1\subseteq E_2\subseteq\cdots\subseteq X$, then $\Hh^h\bigl(\bigcup_nE_n\bigr)=\lim_{n\to\infty}\Hh^h(E_n)$.
\end{corollary}

\begin{proof}
By \cref{thm:borelregular} choose Borel sets $B_k\supseteq E_k$ with $\Hh^h(B_k)=\Hh^h(E_k)$, and put
$C_n=\bigcap_{k\ge n}B_k$. Each $C_n$ is Borel, $C_n\subseteq C_{n+1}$, and $E_n\subseteq C_n\subseteq B_n$
(because $E_n\subseteq E_k\subseteq B_k$ for $k\ge n$); hence $\Hh^h(C_n)=\Hh^h(E_n)$. Since $\Hh^h$ is
countably additive on Borel sets (\cref{thm:metric}), it is continuous along increasing sequences of
Borel sets, and therefore
\[
\Hh^h\Bigl(\bigcup_nE_n\Bigr)\le\Hh^h\Bigl(\bigcup_nC_n\Bigr)=\lim_n\Hh^h(C_n)=\lim_n\Hh^h(E_n)
\le\Hh^h\Bigl(\bigcup_nE_n\Bigr),
\]
where the limits exist by monotonicity.
\end{proof}

Borel regularity must not be confused with approximation from above by \emph{open} sets.

\begin{example}[No outer regularity by open sets]\label{ex:noopen}
Let $X=\R$ and $h(x,r)=r^{1/2}$. Then $\Hh^h(\{0\})=0$ by \cref{rem:diamzero}, but $\Hh^h(U)=\infty$ for
every open set $U\ni0$. Indeed, $U$ contains a nondegenerate compact interval $I$. Let $0<\delta<1$ and let
$(A_j)$ be a $\delta$-covering of $I$, with $d_j=\diam A_j\le\delta$. Since $d_j^{1/2}\ge\delta^{-1/2}d_j$ for
$0\le d_j\le\delta$, and by the lower bound in the proof of \cref{lem:lebesgue},
\[
\sum_jc_h(A_j)=\sum_jd_j^{1/2}\ge\delta^{-1/2}\sum_jd_j\ge\delta^{-1/2}\,\lambda^*(I)=\delta^{-1/2}|I| .
\]
Hence $\Hh^h_\delta(I)\ge\delta^{-1/2}|I|\to\infty$ as $\delta\downarrow0$, and $\Hh^h(U)\ge\Hh^h(I)=\infty$.
\end{example}

\begin{remark}
For spatial gauges, approximation by open sets is delicate for a structural reason as well: replacing a
covering set $A$ by an open neighbourhood $A_\rho=\{x:\dist(\{x\},A)<\rho\}$ increases the diameter (up to
$\diam A+2\rho$), and since $h$ need not be monotone or continuous in $r$, the cost $c_h(A_\rho)$ is not
controlled by $c_h(A)$. Continuity in the spatial variable and suitable uniform control in the radius are
natural hypotheses for results of this type; we do not pursue them here.
\end{remark}

\section{The local stability problem}\label{sec:local}

The uniform comparison in \cref{thm:uniform} is strong. It is natural to ask whether it can be weakened to
pointwise asymptotic equivalence.

\begin{problem}[Pointwise-to-measure stability]\label{prob:pointwise}
Suppose $h,g\colon X\times(0,r_0)\to(0,\infty)$ satisfy
\[
\lim_{r\downarrow0}\frac{g(x,r)}{h(x,r)}=1\qquad\text{for every fixed }x\in X.
\]
Find additional hypotheses, as weak as possible, under which $\Hh^g(E)=\Hh^h(E)$ for every Borel set
$E\subseteq X$.
\end{problem}

The difficulty is that the covering cost is an infimum over \emph{all} points of a covering set. Pointwise
convergence at each fixed point does not provide a comparison uniform over all points that may occur in
a small covering set. The following example shows that additional hypotheses are indeed indispensable,
even on a compact space and for bounded gauges.

\begin{proposition}[Pointwise equivalence does not suffice]\label{prop:counterexample}
Let $X=[0,1]$ with the Euclidean metric and fix $\kappa\in(0,1)$. For $0<r<1$ let $D_r$ be the set of
rational numbers $p/q\in[0,1]$, written in lowest terms with $q\ge1$, such that $q\ge1/r$. Define
\[
h(x,r)=r,\qquad
g(x,r)=\begin{cases}\kappa\,r,& x\in D_r,\\ r,& x\notin D_r,\end{cases}\qquad(x\in[0,1],\ 0<r<1).
\]
Then $0<g\le h\le1$, and for every $x\in[0,1]$ there is $\rho(x)>0$ with $g(x,r)=h(x,r)$ for all
$0<r<\rho(x)$; in particular $g(x,r)/h(x,r)\to1$ for every $x$. Nevertheless
\[
\Hh^g(E)=\kappa\,\lambda^*(E)=\kappa\,\Hh^h(E)\qquad\text{for every }E\subseteq[0,1];
\]
in particular $\Hh^g([0,1])=\kappa<1=\Hh^h([0,1])$.
\end{proposition}

\begin{proof}
\emph{Pointwise behaviour.} If $x$ is irrational, then $x\notin D_r$ for all $r$, and we take
$\rho(x)=1$. If $x=p/q$ in lowest terms, then $x\in D_r$ if and only if $r\ge1/q$; hence $g(x,r)=h(x,r)$
for $0<r<\rho(x):=1/q$.

\emph{A density property.} Let $P\subseteq[0,1]$ be an interval of length $r\in(0,1)$. Then $P$ contains a
point of $D_r$. Indeed, $P$ contains infinitely many rational numbers, whereas $[0,1]$ contains only
finitely many rationals whose reduced denominator is $<1/r$ (at most $q+1$ for each of the finitely
many $q<1/r$). Consequently $c_g(P)=\inf_{x\in P}g(x,r)=\kappa r=\kappa\,c_h(P)$, because $g\ge\kappa h$
everywhere and the value $\kappa r$ is attained on $P\cap D_r\neq\emptyset$.

\emph{Lower bound.} Since $\kappa h\le g\le h$ on $X\times(0,1)$, \cref{cor:comparable} gives
$\kappa\Hh^h(E)\le\Hh^g(E)\le\Hh^h(E)$, and $\Hh^h=\lambda^*$ on $[0,1]$ by \cref{lem:lebesgue}.

\emph{Upper bound.} Let $E\subseteq[0,1]$, $0<\delta<1$ and $\eps>0$. As in the proof of
\cref{lem:lebesgue}, there are closed intervals $P_{k,i}$ of length at most $\delta$ such that the sets
$Q_{k,i}=P_{k,i}\cap[0,1]$ cover $E$ and $\sum_{k,i}|P_{k,i}|\le\lambda^*(E)+\eps$. Each $Q_{k,i}$ is an
interval contained in $[0,1]$; if it is degenerate, its cost is $0$, and otherwise its length
$r=\diam Q_{k,i}\in(0,\delta]$ and, by the density property, $c_g(Q_{k,i})=\kappa\,\diam Q_{k,i}$. Hence
\[
\Hh^g_\delta(E)\le\sum_{k,i}c_g(Q_{k,i})\le\kappa\sum_{k,i}|P_{k,i}|\le\kappa(\lambda^*(E)+\eps).
\]
Letting $\eps\downarrow0$ and then $\delta\downarrow0$ gives $\Hh^g(E)\le\kappa\lambda^*(E)$, which together
with the lower bound proves the claim.
\end{proof}

In this example the ratio $g(\cdot,r)/h(\cdot,r)$ is discontinuous: for every $r$ it equals $\kappa$ on the
dense set $D_r$. For ratios that are continuous in the spatial variable, a Baire-category argument
converts pointwise convergence into locally uniform convergence on a large set.

\begin{proposition}[Baire localization]\label{prop:baire}
Let $X$ be a Baire space (for instance a complete metric space), and let
$h,g\colon X\times(0,r_0)\to(0,\infty)$ be such that, for each $r\in(0,r_0)$, the ratio
$\varphi_r=g(\cdot,r)/h(\cdot,r)$ is continuous on $X$, and $\varphi_r(x)\to1$ as $r\downarrow0$ for every
$x\in X$. For $\eta>0$ let $G_\eta$ be the set of points $x\in X$ having an open neighbourhood $U$ and a
number $\rho\in(0,r_0)$ such that $|\varphi_r-1|\le\eta$ on $U$ for all $0<r<\rho$. Then $G_\eta$ is open and
dense, and $G=\bigcap_{k\ge1}G_{1/k}$ is a dense $G_\delta$ subset of $X$.
\end{proposition}

\begin{proof}
$G_\eta$ is open, because the neighbourhood $U$ in its definition is contained in $G_\eta$. To prove
density, let $V\neq\emptyset$ be open. For $n\in\N$ with $1/n<r_0$ put
\[
F_n=\bigl\{x\in V:\ |\varphi_r(x)-1|\le\eta\ \text{for all } 0<r<1/n\bigr\}
=V\cap\bigcap_{0<r<1/n}\varphi_r^{-1}\bigl([1-\eta,1+\eta]\bigr).
\]
By continuity of each $\varphi_r$, every $F_n$ is relatively closed in $V$; by pointwise convergence,
$V=\bigcup_nF_n$. A nonempty open subspace of a Baire space is a Baire space, hence is not a countable
union of nowhere dense subsets; therefore some $F_n$ has nonempty interior $W$ relative to $V$, and $W$ is
open in $X$ because $V$ is. Every point of $W$ belongs to $G_\eta$ (take $U=W$ and $\rho=1/n$), so
$G_\eta\cap V\ne\emptyset$. Finally, $G$ is a countable intersection of dense open sets, hence a dense
$G_\delta$ set by the Baire property.
\end{proof}

\begin{theorem}[Local stability on compact sets]\label{thm:compactlocal}
Under the hypotheses of \cref{prop:baire}, $\Hh^g(K)=\Hh^h(K)$ for every compact set $K\subseteq G$, and
more generally $(1-\eta)\Hh^h(K)\le\Hh^g(K)\le(1+\eta)\Hh^h(K)$ for every compact $K\subseteq G_\eta$ with
$0<\eta<1$. Consequently $\Hh^g(E)=\Hh^h(E)$ for every set $E$ that is a countable increasing union of compact
subsets of $G$.
\end{theorem}

\begin{proof}
Let $0<\eta<1$ and let $K\subseteq G_\eta$ be compact and nonempty. Each $x\in K$ has an open
neighbourhood $U_x$ and $\rho_x>0$ as in the definition of $G_\eta$; by compactness,
$K\subseteq U:=U_{x_1}\cup\dots\cup U_{x_N}$, and with $\rho=\min_i\rho_{x_i}$ we have
$|\varphi_r-1|\le\eta$ on $U$ for $0<r<\rho$. Let $\tau=\dist(K,X\setminus U)$ if $U\ne X$ and $\tau=\infty$
otherwise; $\tau>0$ because $K$ is compact, $X\setminus U$ is closed and disjoint from $K$, and the
continuous function $\dist(\cdot,X\setminus U)$ attains a positive minimum on $K$.

Fix $0<\delta<\min\{\rho,\tau\}$ and let $(A_j)$ be a $\delta$-covering of $K$. Discard the sets that do not
meet $K$; the remaining family $(A_j)_{j\in J}$ still covers $K$ and has smaller total cost. Each $A_j$,
$j\in J$, is contained in $U$: if $x\in A_j\cap K$ and $y\in A_j$, then $d(x,y)\le\delta<\tau$, so
$y\in U$. Hence for $j\in J$ with $r=\diam A_j>0$ we have $(1-\eta)h(y,r)\le g(y,r)\le(1+\eta)h(y,r)$ for all
$y\in A_j$, and taking infima over $y\in A_j$ as in the proof of \cref{lem:fixedscale},
$(1-\eta)c_h(A_j)\le c_g(A_j)\le(1+\eta)c_h(A_j)$ for all $j\in J$. Therefore
\[
\sum_jc_g(A_j)\ge\sum_{j\in J}c_g(A_j)\ge(1-\eta)\sum_{j\in J}c_h(A_j)\ge(1-\eta)\Hh^h_\delta(K),
\]
\[
\Hh^g_\delta(K)\le\sum_{j\in J}c_g(A_j)\le(1+\eta)\sum_{j\in J}c_h(A_j)\le(1+\eta)\sum_jc_h(A_j).
\]
Taking infima over all $\delta$-coverings of $K$ gives
$(1-\eta)\Hh^h_\delta(K)\le\Hh^g_\delta(K)\le(1+\eta)\Hh^h_\delta(K)$ for all small $\delta$, and letting
$\delta\downarrow0$ yields the two-sided estimate. If $K\subseteq G$, the estimate holds for every
$\eta=1/k$, and letting $k\to\infty$ gives $\Hh^g(K)=\Hh^h(K)$ (if one side is infinite, so is the other).
The last assertion follows from \cref{cor:increasing}, applied to $\Hh^g$ and to $\Hh^h$.
\end{proof}

\Cref{thm:compactlocal} pinpoints where the difficulty of \cref{prob:pointwise} lies for continuous ratios:
in the closed nowhere dense set $X\setminus G_\eta$. Such a set may nevertheless carry positive
$\Hh^h$-measure (think of a Cantor set of positive Lebesgue measure in $[0,1]$), so a positive answer
requires controlling coverings near $X\setminus G_\eta$.

\begin{problem}[Optimality and counterexamples]\label{prob:optimal}
By \cref{prop:counterexample}, pointwise asymptotic equivalence does not imply equality of the induced
outer measures, even when $X$ is compact and the gauges are uniformly bounded. Decide whether the
implication holds when, in addition, the ratio $g(\cdot,r)/h(\cdot,r)$ is continuous on $X$ for each $r$
(or $g$ and $h$ are jointly continuous), and, if it fails, which further hypotheses restore it.
\end{problem}

A successful solution should contain both a positive theorem and a counterexample showing that its
hypotheses cannot be substantially weakened. Several approaches deserve investigation:
\begin{enumerate}[label=(\roman*)]
\item \emph{Compact-uniform convergence.} Assume that the ratio converges uniformly in $x$ on every compact
subset. The proof of \cref{thm:compactlocal} shows that the points of a covering set lying outside a
neighbourhood of the compact set never need to be considered at small scales; it then remains to pass from
compact sets to Borel sets, for instance under inner regularity of $\Hh^h$ by compact sets.
\item \emph{Equicontinuity in the spatial variable.} A common modulus of continuity for the ratios
$\varphi_r$ converts pointwise convergence into uniform convergence on compact sets, by a finite-net
argument, and reduces the problem to (i). Without equicontinuity, \cref{prop:baire} still yields locally
uniform convergence on the dense open sets $G_\eta$.
\item \emph{Restricted covering families.} If the construction is modified so that the location parameter
is selected from the portion of the covering set lying in a prescribed compact region, localization may
become possible. This defines a different construction and must be studied on its own terms.
\item \emph{Capacity or tightness assumptions.} A condition preventing low-cost points from being used by
coverings of geometrically distant regions may yield a useful substitute for uniform convergence; note that
in \cref{prop:counterexample} the low-cost sets $D_r$ are dense at every scale.
\end{enumerate}

\section{Relation to previous work and scope of the results}\label{sec:literature}

The construction of metric outer measures by coverings is classical. Carathéodory's method,
the metric outer-measure principle, and the criterion ensuring measurability of closed sets are
treated in standard references such as \cite[\S\S~2.3, 2.10]{Federer1969},
\cite{Rogers1970}, \cite[Chapters~1 and~4]{Mattila1995}, \cite{Folland1999},
\cite[Chapter~5]{Taylor2006}, and \cite{EoMOuter}. Accordingly, the metric property proved
in \cref{thm:metric} is included to make the framework self-contained; it is not asserted as a
new general principle. The same qualification applies to the elementary comparison argument
in \cref{thm:uniform}. The closure inequality in \cref{lem:closed}, which follows from taking
the location parameter over the covering set itself, is the simple structural fact that makes
the Borel-hull argument particularly direct in this setting.

Spatially dependent Hausdorff gauges have appeared in several forms. Gauges of the type
$h(x,r)=r^t\nu(B(x,r))^q$, where $\nu$ is a measure, occur in Olsen's multifractal formalism
\cite{Olsen1995}. Centred Hausdorff measures, introduced by Saint Raymond and Tricot
\cite{SaintRaymondTricot1988}, were studied further by Schechter
\cite{Schechter2000,Schechter2002}. In centred constructions the admissible centre of a
covering ball is constrained to belong to the set being measured. That constraint makes
replacement of a covering set by its closure more delicate. In the construction used here,
the location variable ranges over the covering set, and hence closure does not increase the
cost. This distinction explains the unconditional Borel-regularity statement in
\cref{thm:borelregular}; it should not be confused with a corresponding assertion for
centred Hausdorff measures.

The stability questions considered here are different from the basic construction problem.
The logarithmic criterion in \cref{thm:log} is a quantitative consequence of uniform
comparison, and its proof is elementary. Its role is to identify the exact scale at which
uniform exponent perturbations produce a nontrivial multiplicative change. The more
substantial issue is the failure of pointwise stability: \cref{prop:counterexample} shows
that even eventual pointwise equality of the gauges at each fixed spatial point does not
prevent a fixed change in the induced measure. The mechanism is that every interval at a
given scale contains points where the covering cost is reduced, although the low-cost
points depend on the scale.

The continuity assumption in \cref{prop:baire} prevents that particular dense-set
mechanism from operating without restriction. Baire's theorem gives local uniform control
on a dense open set at each fixed tolerance, and compactness then turns this control into
comparison of covering costs. The resulting statement, \cref{thm:compactlocal}, applies
to compact subsets of a dense $G_\delta$ set, not to all Borel sets. This limitation is
essential to the formulation: a meagre exceptional set can still have positive measure,
and the present argument does not control coverings that approach that exceptional set.

Thus the contribution of the paper is best understood as a set of precise stability
statements and a counterexample delimiting their scope, rather than as a new construction
of metric outer measure. The following question remains open.

\section{Conclusions}

We have studied the stability of Carathéodory outer measures generated by covering costs
whose values depend on both scale and spatial position. The underlying construction yields
a metric outer measure and is Borel regular; these foundational facts are included to clarify
the setting in which the stability questions are posed. In particular, every set has an
$F_{\sigma\delta}$ Borel hull of the same outer measure, whereas approximation from above
by open sets may fail.

For uniformly controlled perturbations, the induced outer measures satisfy sharp comparison
principles. In the power-gauge setting, the condition
\[
\sup_{\substack{x\in X\\0<r\le\delta}}
\left|a(x,r)\log\frac1r-\lambda\right|\longrightarrow0
\]
implies $\mathcal H^g=e^{-\lambda}\mathcal H^h$. The logarithmic scale therefore identifies
the threshold for this uniform criterion.

The principal nonuniform phenomenon is that pointwise asymptotic equivalence does not imply
equality of the induced outer measures. The example on $[0,1]$ gives bounded gauges whose
ratio is eventually one at each fixed point, while the resulting measures differ by a fixed
factor. Under spatial continuity of the ratio, Baire-category arguments recover equality on
compact subsets of a dense $G_\delta$ set. They do not, however, settle the question for all
Borel sets, because the meagre exceptional set may carry positive measure.

A natural next step is to identify additional hypotheses—such as equicontinuity, compact-uniform
convergence, or suitable control of the exceptional set—that yield stability on larger classes
of sets. Any such result should be accompanied by counterexamples that clarify the sharpness
of its assumptions.


\end{document}